\documentclass[a4paper;11pt]{amsart}
\usepackage{mathrsfs}\usepackage{graphicx}
 \usepackage[arrow,matrix]{xy}
\usepackage{amsmath,amssymb,amscd,bbm,amsthm,mathrsfs}
\usepackage{color,xcolor}
\usepackage{graphicx}
\usepackage{manfnt}
\newtheorem{thm}{Theorem}[section]
\newtheorem{lem}{Lemma}[section]
\newtheorem{cor}{Corollary}[section]

\newtheorem{rem}{Remark}[section]

\theoremstyle{definition}

\begin{document}
\numberwithin{equation}{section}

\title[]{On Pinned Falconer Distance Problem
for Cartesian Product Sets: the Parabolic Method}

\author{Ji Li, Chong-Wei Liang and Chun-Yen Shen}

\maketitle

\begin{center}
\begin{minipage}{13.5cm}\small

{\noindent  {\bf Abstract:}\  
 The Falconer distance problem 
for Cartesian product sets was introduced and studied by Iosevich and Liu (\cite{MR3525385}).
  In this
  paper, by implementing a new observation on Cartesian product sets
associated with a particular parabolic structure, we study the pinned version of Falconer distance problem 
for Cartesian product sets, and improve the threshold for the Falconer distance set in \cite{MR3525385} in certain case. 
}

\end{minipage}
\end{center}


\section{Introduction and Statement of Main Results}
The Falconer distance conjecture (\cite{MR0834490}) says that if the Hausdorff dimension of $E \subset {\Bbb R}^d$, $d \ge 2$, is greater than $\frac{d}{2}$, then the Lebesgue measure of the distance set $\Delta(E)=\{|x-y|: x,y \in E \}$ is positive. 

Recent celebrated results \cite{MR4055179,MR4297185,du2024newimprovementfalconerdistance} show that for every 
  compact subset $E$ of ${\Bbb R}^d$ with $d\geq2$,  the Lebesgue measure of the
  distance set $\Delta(E)=\{|x-y|: x,y \in E \}$ is positive if the
  Hausdorff dimension of $E$ satisfies
  ${\rm dim}_{\mathcal H}(E)>{d\over2}+{1\over4}$ when $d$ is even, and ${\rm dim}_{\mathcal H}(E)>{d\over2}+{1\over4}+{1\over 8d-4}$ when $d$ is odd. 

  This improved the well-known result by Wolff \cite{MR1692851} in two dimensions and Erdogan \cite{MR2152236} in higher
dimensions.

Recall that Falconer distance problem on Cartesian product sets 
was studied by Iosevich and Liu (\cite{MR3525385}) via Mattila integral,
which states as follows.

\vskip.1cm
\noindent {\bf Theorem A} (\cite{MR3525385}){\bf.}\  {\it Let $E=A \times B$, where $A$ and $B$
  are compact subsets of ${\Bbb R}$ with positive $s_A, s_B$-dimensional Hausdorff measure, respectively. If $s_A+s_B+\max(s_A,s_B)>2$, the Lebesgue measure of
  $\Delta(E)$ is positive. 
}
\vskip.1cm


\noindent {\bf Theorem B} (\cite{MR3525385}){\bf.}\  {\it  Suppose that $E$ is a compact subset of ${\Bbb R}^d$ of the form $A_1 \times A_2 \times \dots \times A_d$, where $A_j \subset {\Bbb R}$ has positive $s_j$-dimensional Hausdorff measure for all $1\leq j\leq d$. Suppose that $\sum_{j=1}^d s_j>\frac{d^2}{2d-1}$. Then the Lebesgue measure of $\Delta(E)$ is positive. 
}
\vskip.1cm

Their result improved Erdogan's $\frac{d}{2}+\frac{1}{3}$ exponent in higher dimension for Cartesian products.  Note also that they studied the pinned version of Falconer distance problem 
 \cite{MR3917228}.

\vskip.1cm
In this paper, we study the pinned version of Falconer distance problem 
for Cartesian product sets
by implementing a new observation on Cartesian product sets
associated with a  particular parabolic structure. 
We improve the threshold for the Falconer distance set of product sets $A_1 \times A_2 \times \dots \times A_d$ in Theorem {\bf B} in certain case.


We now state our main results in details.

In what follows,  for any set $E\subset\mathbb R^d$ and let $x\in E$, we define $\triangle_x(E)$ to be the pinned version of the distance set, that is
 \begin{align*}
     \triangle_x(E):=\left\{|x-y|:y\in E\right\}.
 \end{align*} 
Moreover, for a set $A \subset \mathbb R$, we denote $A^d= A \times\cdots\times A\subset\mathbb R^d$.

\begin{thm}\label{main1} Let $A, B\subset\mathbb R$ be compact subsets and $d\geq3$.
\noindent 1. If one of  the following conditions holds
 \begin{align*}
\begin{cases}
\displaystyle     1<\dim_{\mathcal{H}}((A\cap B)^2)\leq\frac{57}{56}\\[9pt]
\displaystyle     \dim_{\mathcal{H}}(A)+\frac{(d-2)}{d}\dim_{\mathcal{H}}(B)\geq\frac{21d^2-71d+21\beta\cdot d+29}{21d(d-2)}
\end{cases}
{\rm\ \  or\ \ }
\end{align*}
\begin{align*}
\begin{cases}
\displaystyle
    \frac{57}{56}<\dim_{\mathcal{H}}((A\cap B)^2)\leq\frac{5}{4}\\[9pt]
\displaystyle   \left[\frac{(d-2)}{d}\dim_{\mathcal{H}}(B)+\dim_{\mathcal{H}}(A)\right]+\frac{16d-16}{3d(d-2)}\cdot\dim_{\mathcal{H}}(A\cap B)> 1+\frac{3d\beta+4d-4}{3d(d-2)}
\end{cases}
{\rm\ \  or\ \ }
\end{align*}
\begin{align*}
\begin{cases}
\displaystyle
   \frac{5}{4}<\dim_{\mathcal{H}}((A\cap B)^2)\\[9pt]
\displaystyle \dim_{\mathcal{H}}(A)+\frac{(d-2)}{d}\dim_{\mathcal{H}}(B)>1+\frac{\beta d-2d+2}{d(d-2)},  
\end{cases}
\end{align*}
then there is a point $(b_1,\ldots,b_d)\in B^d$ such that the Hausdorff dimension of pinned distance set of the product $\triangle_{(b_1,\ldots,b_d)}(A^d)$ is no less than $\beta$, where $A^d$ is the product set of $A$.

\bigskip

\noindent 2. 
If one of the following conditions holds
\begin{align*}
\begin{cases}
\displaystyle     1<\dim_{\mathcal{H}}((A\cap B)^2)\leq\frac{57}{56}\\[9pt]
\displaystyle     \dim_{\mathcal{H}}(A)+\frac{(d-2)}{d}\dim_{\mathcal{H}}(B)\geq\frac{21d^2-29d+29}{21d(d-2)}
\end{cases}
{\rm\ \  or\ \ }
\end{align*}
\begin{align*}
\begin{cases}
\displaystyle
    \frac{57}{56}<\dim_{\mathcal{H}}((A\cap B)^2)\leq\frac{5}{4}\\[9pt]
\displaystyle   \left[\frac{(d-2)}{d}\dim_{\mathcal{H}}(B)+\dim_{\mathcal{H}}(A)\right]+\frac{16d-16}{3d(d-2)}\cdot\dim_{\mathcal{H}}(A\cap B)> 1+\frac{10d-4}{3d(d-2)}
\end{cases}
{\rm\ \  or\ \ }
\end{align*}
\begin{align*}
\begin{cases}
\displaystyle
   \frac{5}{4}<\dim_{\mathcal{H}}((A\cap B)^2)\\[9pt]
\displaystyle \dim_{\mathcal{H}}(A)+\frac{(d-2)}{d}\dim_{\mathcal{H}}(B)>1+\frac{2}{d(d-2)},  
\end{cases}
\end{align*}
then there is $(b_1,\ldots,b_d)\in B^d$ such that  $\triangle_{(b_1,\ldots,b_d)}(A^d)$ contains an non-empty interval, where $A^d$ is the product set of $A$.

\bigskip

\noindent 3. 
If one of the following conditions holds\begin{align*}
\begin{cases}
\displaystyle     1<\dim_{\mathcal{H}}((A\cap B)^2)\leq\frac{57}{56}\\[9pt]
\displaystyle     \dim_{\mathcal{H}}(A)+\frac{(d-2)}{d}\dim_{\mathcal{H}}(B)\geq\frac{(21d-29)(d-1)}{21d(d-2)}
\end{cases}
{\rm\ \  or\ \ }
\end{align*}
\begin{align*}
\begin{cases}
\displaystyle
    \frac{57}{56}<\dim_{\mathcal{H}}((A\cap B)^2)\leq\frac{5}{4}\\[9pt]
\displaystyle   \left[\frac{(d-2)}{d}\dim_{\mathcal{H}}(B)+\dim_{\mathcal{H}}(A)\right]+\frac{16d-16}{3d(d-2)}\cdot\dim_{\mathcal{H}}(A\cap B)> 1+\frac{7d-4}{3d(d-2)}
\end{cases}
{\rm\ \  or\ \ }
\end{align*}
\begin{align*}
\begin{cases}
\displaystyle
   \frac{5}{4}<\dim_{\mathcal{H}}((A\cap B)^2)\\[9pt]
\displaystyle \dim_{\mathcal{H}}(A)+\frac{(d-2)}{d}\dim_{\mathcal{H}}(B)>1-\frac{1}{d},  
\end{cases}
\end{align*}
then there is $(b_1,\ldots,b_d)\in B^d$ such that $\triangle_{(b_1,\ldots,b_d)}(A^d)$ has positive Lebesgue measure, where $A^d$ is the product set of $A$.
\end{thm}
Note that our main theorem implies the following 
special case.

\begin{cor}\label{mainsub1}
     Let $A\subset\mathbb R$ be a compact subset and $d\geq3$.
\noindent 1. If one of  the following conditions holds
 \begin{align*}
\begin{cases}
\displaystyle     1<\dim_{\mathcal{H}}(A^2)\leq\frac{57}{56}\\[9pt]
\displaystyle     \dim_{\mathcal{H}}(A)\geq\frac{21d^2-71d+21\beta\cdot d+29}{42(d-1)(d-2)}
\end{cases}
{\rm\ \  or\ \ }
\begin{cases}
\displaystyle
    \frac{57}{56}<\dim_{\mathcal{H}}(A^2)\leq\frac{5}{4}\\[9pt]
\displaystyle   \dim_{\mathcal{H}}(A)>\frac{3d(d-2+\beta)}{(d-1)(6d+4)}+\frac{2}{3d+2}
\end{cases}
\end{align*}
or
\begin{align*}
\begin{cases}
\displaystyle
   \frac{5}{4}<\dim_{\mathcal{H}}(A^2)\\[9pt]
\displaystyle \dim_{\mathcal{H}}(A)>\frac{d^2-4d+\beta d+2}{(d-2)(2d-2)},  
\end{cases}
\end{align*}
then there is $(y_1,\ldots,y_d)\in A^d$ such that the Hausdorff dimension of pinned distance set of the product $\triangle_{(y_1,\ldots,y_d)}(A^d)$ is no less than $\beta$, where $A^d$ is the product set of $A$.

\bigskip

\noindent 2. 
If one of the following conditions holds
\begin{align*}
\begin{cases}
\displaystyle
     1<\dim_{\mathcal{H}}(A^2)\leq\frac{57}{56}\\[9pt]
 \displaystyle    \dim_{\mathcal{H}}(A)\geq\frac{1}{2}+\frac{34d-13}{42(d-1)(d-2)}
\end{cases}
,\text{or}\,\,
\begin{cases}
\displaystyle    \frac{57}{56}<\dim_{\mathcal{H}}(A^2)\leq\frac{5}{4}\\[9pt]
\displaystyle   \dim_{\mathcal{H}}(A)>\frac{1}{2}+\frac{5d-2}{(d-1)(6d+4)}
\end{cases}
\end{align*}
or
\begin{align*}
\begin{cases}
\displaystyle   \frac{5}{4}<\dim_{\mathcal{H}}(A^2)\\[9pt]
\displaystyle \dim_{\mathcal{H}}(A)>\frac{1}{2}+\frac{d}{(d-2)(2d-2)},  
\end{cases}
\end{align*}
then there is $(y_1,\ldots,y_d)\in A^d$ such that  $\triangle_{(y_1,\ldots,y_d)}(A^d)$ contains an non-empty interval, where $A^d$ is the product set of $A$.

\bigskip

\noindent 3. 
If one of the following conditions holds
\begin{align*}
\begin{cases}
\displaystyle     1<\dim_{\mathcal{H}}(A^2)\leq\frac{57}{56}\\[9pt] 
\displaystyle     \dim_{\mathcal{H}}(A)\geq\frac{1}{2}+\frac{13}{42(d-2)}
\end{cases}
,\text{or}\,\,
\begin{cases}
\displaystyle    \frac{57}{56}<\dim_{\mathcal{H}}(A^2)\leq\frac{5}{4}\\[9pt]
\displaystyle   \dim_{\mathcal{H}}(A)>\frac{1}{2}+\frac{1}{3d+2}
\end{cases}
\end{align*}
or
\begin{align*}
\begin{cases}
  \displaystyle \frac{5}{4}<\dim_{\mathcal{H}}(A^2)\\[9pt]
\displaystyle \dim_{\mathcal{H}}(A)>\frac{1}{2},  
\end{cases}
\end{align*}
then there is $(y_1,\ldots,y_d)\in A^d$ such that $\triangle_{(y_1,\ldots,y_d)}(A^d)$ has positive Lebesgue measure, where $A^d$ is the product set of $A$.

\end{cor}

Based on part 3 in the Corollary \ref{mainsub1}  on the pinned version of the distance set, we improve the threshold for the Falconer distance set of product sets $A_1 \times A_2 \times \dots \times A_d$ in \cite{MR3525385} in certain case: the Hausdorff dimensions $s_{d-1}$ and $s_d$
of $A_{d-1}$ and $A_d$, respectively, satisfying 
$s_{d-1}+s_d>\frac{5}{4}$.
We state this in details in below.
\begin{cor}\label{sec}
    Let $d\geq3$. Suppose that $E$ is a compact subset of ${\Bbb R}^d$ of the form $A_1 \times A_2 \times \dots \times A_d$, where $A_j \subset {\Bbb R}$ has positive $s_j$-dimensional Hausdorff measure for all $1\leq j\leq d$. Suppose that $\sum_{j=1}^{d-2} s_j>\frac{d}{2}-1$ and $s_{d-1}+s_d>\frac{5}{4}$. Then there is a point $(y_1,\ldots,y_d)\in E$ such that the Lebesgue measure of the pinned distance set $\Delta_{(y_1,\ldots,y_d)}(E)$ is positive.
\end{cor}

\begin{rem}
    Note that Iosevich and Liu (\cite{MR3917228})
    studied the pinned version of distance set
    for general set $E,F\subset \mathbb R^d$.
    Our results improve the threshold obtained by \cite{MR3917228} when applying their results to $E=A_1\times\cdots\times A_d$ and $F=B_1\times\cdots\times B_d$. In fact, one can see from our proofs below that as long as one can prove a variant of distance result for the particular parabolic distance, one can improve the result for the classical distance for product sets. 
    
    \end{rem}

Besides, there are other variant distance problems related to the classical Falconer distance set. (see for example \cite{MR3365800, MR4289240, MR4439497, MR4490711, MR4623515, MR4530888, borges2023singularvariantfalconerdistance, MR4846807}) Our technique can be also applied to improve the threshold for those distance problems in the case of Cartesian product sets.

\section{Proofs of Theorem \ref{main1} and Corollary \ref{mainsub1}: Parabolic Attack}
Let $\Phi:\mathbb R^d\times\mathbb R^d\longrightarrow\mathbb R$ be a smooth function that satisfies \textbf{Phong--Stein rotation curvature condition} and \textbf{Sogge's cinematic curvature condition}, that is 
$\Phi$ has a nonzero Monge--Ampere determinant
\begin{align*}
    \det \begin{pmatrix}
0 &\quad \nabla_x\Phi \\\nabla_y\Phi&\quad
\frac{\partial^2\Phi}{\partial x\partial y}\\[5pt]
\end{pmatrix}\neq0
\end{align*}
and for any $t>0$, $x\in\mathbb R^d$, $\left\{\nabla_y\Phi: \Phi(x,y)=t\right\}$ has nonzero Gaussian curvature. In particular, the parabolic distance $\Phi(x,y):=(x_1-y_1)^2+\cdots+(x_{d-1}-y_{d-1})^2+(x_d-y_d)$ satisfies both curvature conditions.

\newpage

\begin{thm}[\cite{MR3917228}]\label{thm IL}~\\
    Let $\Phi\in C^\infty(\mathbb R^d\times\mathbb R^d)$, $E,F\subset\mathbb{R}^d$. Suppose that $\Phi$ satisfies the Phong--Stein rotation curvature condition and the cinematic curvature condition. Then there is a probability measure $\mu_F$ on $F$ such that for $\mu_F$-a.e. $y\in F$,
    \begin{align*}
        &(1)\,\,\dim_{\mathcal{H}}(\triangle^\Phi_y(E))\geq\beta,\,\,\text{if}\,\dim_{\mathcal{H}}(E)+\frac{d-1}{d+1}\dim_{\mathcal{H}}(F)>d-1+\beta,\\
          &(2)\,\,\left|\triangle^\Phi_y(E)\right|>0,\,\,\text{if}\,\dim_{\mathcal{H}}(E)+\frac{d-1}{d+1}\dim_{\mathcal{H}}(F)>d,\\
           &(3)\,\,\left(\triangle^\Phi_y(E)\right)^0\neq\emptyset,\,\,\text{if}\,\dim_{\mathcal{H}}(E)+\frac{d-1}{d+1}\dim_{\mathcal{H}}(F)>d+1.
    \end{align*}
\end{thm}
We then have a direct corollary from Theorem \ref{thm IL} as follows when applying their $E$ and $F$ to product sets.
\begin{cor}\label{CoroIL}    
  In particular, if $E=F:=A\times\cdots\times A\subset\mathbb R^d$, where $A\subset\mathbb R$ is a compact subset. Then there is a probability measure $\mu_F$ on $F$ such that for $\mu_F$-a.e $y\in F$,
     \begin{align*}
    &(1)\,\,\dim_{\mathcal{H}}(A)>\frac{(d-1+\beta)(d+1)}{2d^2}
    \implies \,\dim_{\mathcal{H}}(\triangle^\Phi_y(E))\geq\beta,\\
    &(2)\,\,\dim_{\mathcal{H}}(A)>\frac{d+1}{2d}
    \implies \,\left|\triangle^\Phi_y(E)\right|>0,\\
    &(3)\,\,\dim_{\mathcal{H}}(A)>\frac{(d+1)^2}{2d^2}
    \implies \,\left(\triangle^\Phi_y(E)\right)^0\neq\emptyset.
    \end{align*}
\end{cor}
 \begin{rem}
    The threshold in Corollary \ref{mainsub1} is better than the threshold in Corollary \ref{CoroIL} for all large $d$.
 \end{rem}
    Next, we recall the following two auxiliary lemmas.
\begin{lem}[\cite{MR4042855}]\label{liu}~\\
For any compact set $\Omega\subset\mathbb R^2$ with $\dim_{\mathcal{H}}(\Omega)>1$, there exists a point $x\in \Omega$ such that
\begin{align*}
 \dim_{\mathcal{H}}(\triangle_x\left(\Omega\right))\geq\min\left\{\frac{4}{3}\cdot\dim_{\mathcal{H}}(\Omega)-\frac{2}{3},\,1\right\}.
\end{align*}
\end{lem}
\begin{lem}[\cite{MR4226184}]\label{sm}~\\
  For any compact set $\Omega\subset\mathbb R^2$ with $\dim_{\mathcal{H}}(\Omega)\in(1,\,1.04)$, there exists a point $x\in \Omega$ such that
\begin{align*}
 \dim_{\mathcal{H}}(\triangle_x\left(\Omega\right))>\frac{29}{42}.
\end{align*}  
\end{lem}
\begin{rem}
    When $\dim_{\mathcal{H}}(\Omega)$ is closed to $1$, then the lower bound $\frac{29}{42}$ is bigger than the lower bound $\min\left\{\frac{4}{3}\dim_{\mathcal{H}}(\Omega)-\frac{2}{3},\,1\right\}$. To compare these two lemmas, one can actually have that there is a point $x\in\Omega$ such that
    \begin{align}
     &\dim_{\mathcal{H}}(\triangle_x\left(\Omega\right))>\frac{29}{42},\,\,\text{if}\,\,\, 1<\dim_{\mathcal{H}}(\Omega)\leq\frac{57}{56};\,\text{and}\notag\\[9pt] 
     &\dim_{\mathcal{H}}(\triangle_x\left(\Omega\right))\geq\begin{cases}\displaystyle
     \frac{4}{3}\dim_{\mathcal{H}}(\Omega)-\frac{2}{3},\,\,\text{if}\,\,\,\frac{57}{56}<\dim_{\mathcal{H}}(\Omega)\leq\frac{5}{4},\\[9pt]
     \displaystyle 1,\,\,\text{if}\,\,\,\frac{5}{4}<\dim_{\mathcal{H}}(\Omega).
     \end{cases}
    \end{align}
\end{rem}

\subsection{Proof of Theorem \ref{main1}}
 Let $d\geq3$ and $y_{d-1}, y_d\in A\cap B$ and $x_0, x_1\in B$. Consider the sets in $\mathbb{R}^{d-1}$.
    \begin{align*}
        &E:=\left(A\times\cdots \times A\right)\times \triangle^2_{(y_{d-1},y_{d})}(A^2)\subset\mathbb R^{d-2}\times\mathbb R^1;\,\text{and}\\
        &F:=\left(B\times\cdots \times B\right)\times -\triangle^2_{(x_0,x_1)}(B^2)\subset\mathbb R^{d-2}\times\mathbb R^1.
    \end{align*}
    If $\Phi_{d-1}:\mathbb R^{d-1}\times\mathbb R^{d-1}\longrightarrow\mathbb R$ is the parabolic distance,  then there is a probability measure $\mu_F$ on $F$ such that for $\mu_F$-a.e. $y\in F$ which is of the form  $$\left(b_1,\ldots,b_{d-2},-\left(|x_0-b_{d-1}|^2+|x_1-b_d|^2\right)\right),$$
    satisfying
     \begin{align*}
        &(1)\,\,\dim_{\mathcal{H}}(\triangle^\Phi_y(E))\geq\beta,\,\,\text{if}\,\dim_{\mathcal{H}}(E)+\frac{d-2}{d}\dim_{\mathcal{H}}(F)>d-2+\beta,\\
          &(2)\,\,\left|\triangle^\Phi_y(E)\right|>0,\,\,\text{if}\,\dim_{\mathcal{H}}(E)+\frac{d-2}{d}\dim_{\mathcal{H}}(F)>d-1,\\
           &(3)\,\,\left(\triangle^\Phi_y(E)\right)^0\neq\emptyset,\,\,\text{if}\,\dim_{\mathcal{H}}(E)+\frac{d-2}{d}\dim_{\mathcal{H}}(F)>d.
    \end{align*}
    However note that $$\triangle^\Phi_y(E)=\triangle^2_{(b_1,\ldots,b_{d-2}, y_{d-1}, y_d)}(A^d)+|x_0-b_{d-1}|^2+|x_1-b_d|^2,$$
    which says the set $\triangle^2_{(y_1,\ldots,y_d)}(A^d)$ is a translation of $\triangle^\Phi_y(E)$ so that we have
  \begin{align*}
    &(1)\,\,\dim_{\mathcal{H}}(\triangle^\Phi_y(E))=  \dim_{\mathcal{H}}(\triangle^2_{(b_1,\ldots,b_{d-2}, y_{d-1}, y_d)}(A^d))=\dim_{\mathcal{H}}(\triangle_{(b_1,\ldots,b_{d-2}, y_{d-1}, y_d)}(A^d)),\,\text{and}\\
    &(2)\,\,\left|\triangle^\Phi_y(E)\right|=\left|\triangle^2_{(b_1,\ldots,b_{d-2}, y_{d-1}, y_d)}(A^d)\right|\leq C_A\cdot\left|\triangle_{(b_1,\ldots,b_{d-2}, y_{d-1}, y_d)}(A^d)\right|,\,\text{and}\\
    &(3)\,\,\left(\triangle^\Phi_y(E)\right)^0\neq\emptyset\implies\left(\triangle_{(b_1,\ldots,b_{d-2}, y_{d-1}, y_d)}(A^d)\right)^0\neq\emptyset.
  \end{align*}
Then one has
 \begin{align*}
        &(1)\,\,\dim_{\mathcal{H}}(\triangle_{(b_1,\ldots,b_{d-2}, y_{d-1}, y_d)}(A^d))\geq\beta,\,\,\text{if}\,\dim_{\mathcal{H}}(E)+\frac{d-2}{d}\dim_{\mathcal{H}}(F)>d-2+\beta,\\
          &(2)\,\,\left|\triangle_{(b_1,\ldots,b_{d-2}, y_{d-1}, y_d)}(A^d)\right|>0,\,\,\text{if}\,\dim_{\mathcal{H}}(E)+\frac{d-2}{d}\dim_{\mathcal{H}}(F)>d-1,\\
           &(3)\,\,\left(\triangle_{(b_1,\ldots,b_{d-2}, y_{d-1}, y_d)}(A^d)\right)^0\neq\emptyset,\,\,\text{if}\,\dim_{\mathcal{H}}(E)+\frac{d-2}{d}\dim_{\mathcal{H}}(F)>d.
    \end{align*}
To guarantee the condition $\dim_{\mathcal{H}}(E)+\frac{d-2}{d}\dim_{\mathcal{H}}(F)>d-2+\beta$ holds, we note that
\begin{align}\label{last11}
    &\dim_{\mathcal{H}}(E)+\frac{d-2}{d}\dim_{\mathcal{H}}(F)\notag\\
    &\geq (d-2)\dim_{\mathcal{H}}(A)+\dim_{\mathcal{H}}(\triangle_{(y_{d-1},y_d)}(A^2))\notag\\
    &\qquad+\frac{d-2}{d}\cdot\left((d-2)\dim_{\mathcal{H}}(B)+\dim_{\mathcal{H}}(\triangle_{(x_0,x_1)}(B^2))\right)\notag\\
    &=(d-2)\left[\frac{(d-2)}{d}\dim_{\mathcal{H}}(B)+\dim_{\mathcal{H}}(A)\right]+\dim_{\mathcal{H}}(\triangle_{(y_{d-1},y_d)}(A^2))\notag\\
    &\qquad+\frac{d-2}{d}\dim_{\mathcal{H}}(\triangle_{(x_0,x_1)}(B^2)).
\end{align}
To apply Lemma \ref{liu} and Lemma \ref{sm}, we need to assume that {the Hausdorff dimension of the set $A\cap B$ is at least $\frac{1}{2}$}, and hence we have the following three cases: \\[10pt]
\textbf{Case}$(1)$: Assume $1<\dim_{\mathcal{H}}((A\cap B)^2)\leq\frac{57}{56}$, then (\ref{last11}) reveals that if we choose the point $(y_{d-1},y_d)=(x_0,x_1)\in A\cap B$ such that Lemma \ref{sm} holds for the set $(A\cap B)^2$, then
\begin{align*}
    &\dim_{\mathcal{H}}(E)+\frac{d-2}{d}\dim_{\mathcal{H}}(F)\\
    &\geq (d-2)\left[\frac{(d-2)}{d}\dim_{\mathcal{H}}(B)+\dim_{\mathcal{H}}(A)\right]+\dim_{\mathcal{H}}(\triangle_{(y_{d-1},y_d)}((A\cap B)^2))\notag\\
    &\qquad+\frac{d-2}{d}\dim_{\mathcal{H}}(\triangle_{(x_0,x_1)}((A\cap B)^2))\\
    &> (d-2)\left[\frac{(d-2)}{d}\dim_{\mathcal{H}}(B)+\dim_{\mathcal{H}}(A)\right]+\frac{2d-2}{d}\cdot\frac{29}{42}.   
\end{align*}
Suppose that $
(d-2)\left[\frac{(d-2)}{d}\dim_{\mathcal{H}}(B)+\dim_{\mathcal{H}}(A)\right]+\frac{2d-2}{d}\cdot\frac{29}{42}\geq d-2+\beta$,
then we have 
\begin{align*}
\dim_{\mathcal{H}}(A)+\frac{(d-2)}{d}\dim_{\mathcal{H}}(B)\geq\frac{21d^2-71d+21\beta\cdot d+29}{21d(d-2)},
\end{align*}
whenever $d\geq3$.\\[10pt] 
\textbf{Case}$(2)$: Assume $\frac{57}{56}<\dim_{\mathcal{H}}((A\cap B)^2)\leq\frac{5}{4}$, then (\ref{last11}) reveals that if we choose the point $(y_{d-1},y_d)=(x_0,x_1)$ such that Lemma \ref{liu} holds for the set $(A\cap B)^2$, then
\begin{align*}
   \dim_{\mathcal{H}}(E)+\frac{d-2}{d}\dim_{\mathcal{H}}(F)&\geq (d-2)\left[\frac{(d-2)}{d}\dim_{\mathcal{H}}(B)+\dim_{\mathcal{H}}(A)\right]\\
   &\qquad+\frac{2d-2}{d}\cdot\left(\frac{4}{3}\dim_{\mathcal{H}}((A\cap B)^2)-\frac{2}{3}\right)\\
&\geq(d-2)\left[\frac{(d-2)}{d}\dim_{\mathcal{H}}(B)+\dim_{\mathcal{H}}(A)\right]\\
&\qquad+\frac{16d-16}{3d}\cdot\dim_{\mathcal{H}}(A\cap B)-\frac{4d-4}{3d}.    
\end{align*}
Suppose that $$
(d-2)\left[\frac{(d-2)}{d}\dim_{\mathcal{H}}(B)+\dim_{\mathcal{H}}(A)\right]+\frac{16d-16}{3d}\cdot\dim_{\mathcal{H}}(A\cap B)-\frac{4d-4}{3d}> d-2+\beta,$$
then we have 
\begin{align*}
\left[\frac{(d-2)}{d}\dim_{\mathcal{H}}(B)+\dim_{\mathcal{H}}(A)\right]+\frac{16d-16}{3d(d-2)}\cdot\dim_{\mathcal{H}}(A\cap B)> 1+\frac{3d\beta+4d-4}{3d(d-2)},
\end{align*}
whenever $d\geq3$.\\[10pt] 
\textbf{Case}$(3)$: Assume $\frac{5}{4}<\dim_{\mathcal{H}}(A^2)$, then (\ref{last11}) reveals that if we choose the point $(y_{d-1},y_d)=(x_0,x_1)$ such that Lemma \ref{liu} holds for the set $(A\cap B)^2$, then
\begin{align*}
   \dim_{\mathcal{H}}(E)+\frac{d-2}{d}\dim_{\mathcal{H}}(F)&\geq (d-2)\left[\frac{(d-2)}{d}\dim_{\mathcal{H}}(B)+\dim_{\mathcal{H}}(A)\right]+\frac{2d-2}{d}. 
\end{align*}
Suppose that $$
(d-2)\left[\frac{(d-2)}{d}\dim_{\mathcal{H}}(B)+\dim_{\mathcal{H}}(A)\right]+\frac{2d-2}{d}> d-2+\beta,$$
then we have 
\begin{align*}
\dim_{\mathcal{H}}(A)+\frac{(d-2)}{d}\dim_{\mathcal{H}}(B)>1+\frac{\beta d-2d+2}{d(d-2)},
\end{align*}
whenever $d\geq3$. Therefore, we can conclude that for all $d\geq3$, if
 \begin{align*}
\begin{cases}
\displaystyle     1<\dim_{\mathcal{H}}((A\cap B)^2)\leq\frac{57}{56}\\[9pt]
\displaystyle     \dim_{\mathcal{H}}(A)+\frac{(d-2)}{d}\dim_{\mathcal{H}}(B)\geq\frac{21d^2-71d+21\beta\cdot d+29}{21d(d-2)}
\end{cases}
{\rm\ \  or\ \ }
\end{align*}
\begin{align*}
\begin{cases}
\displaystyle
    \frac{57}{56}<\dim_{\mathcal{H}}((A\cap B)^2)\leq\frac{5}{4}\\[9pt]
\displaystyle   \left[\frac{(d-2)}{d}\dim_{\mathcal{H}}(B)+\dim_{\mathcal{H}}(A)\right]+\frac{16d-16}{3d(d-2)}\cdot\dim_{\mathcal{H}}(A\cap B)> 1+\frac{3d\beta+4d-4}{3d(d-2)}
\end{cases}
{\rm\ \  or\ \ }
\end{align*}
\begin{align*}
\begin{cases}
\displaystyle
   \frac{5}{4}<\dim_{\mathcal{H}}((A\cap B)^2)\\[9pt]
\displaystyle \dim_{\mathcal{H}}(A)+\frac{(d-2)}{d}\dim_{\mathcal{H}}(B)>1+\frac{\beta d-2d+2}{d(d-2)},  
\end{cases}
\end{align*}
then the Hausdorff dimension of pinned distance set of the product $A^d$ is no less than $\beta$, which improve the threshold in the case of product sets. Similarly, we have
for all $d\geq3$, if
\begin{align*}
\begin{cases}
\displaystyle     1<\dim_{\mathcal{H}}((A\cap B)^2)\leq\frac{57}{56}\\[9pt]
\displaystyle     \dim_{\mathcal{H}}(A)+\frac{(d-2)}{d}\dim_{\mathcal{H}}(B)\geq\frac{(21d-29)(d-1)}{21d(d-2)}
\end{cases}
{\rm\ \  or\ \ }
\end{align*}
\begin{align*}
\begin{cases}
\displaystyle
    \frac{57}{56}<\dim_{\mathcal{H}}((A\cap B)^2)\leq\frac{5}{4}\\[9pt]
\displaystyle   \left[\frac{(d-2)}{d}\dim_{\mathcal{H}}(B)+\dim_{\mathcal{H}}(A)\right]+\frac{16d-16}{3d(d-2)}\cdot\dim_{\mathcal{H}}(A\cap B)> 1+\frac{7d-4}{3d(d-2)}
\end{cases}
{\rm\ \  or\ \ }
\end{align*}
\begin{align*}
\begin{cases}
\displaystyle
   \frac{5}{4}<\dim_{\mathcal{H}}((A\cap B)^2)\\[9pt]
\displaystyle \dim_{\mathcal{H}}(A)+\frac{(d-2)}{d}\dim_{\mathcal{H}}(B)>1-\frac{1}{d},  
\end{cases}
\end{align*}
then the  pinned distance set of the product $A^d$ has positive Lebesgue measure; and if
\begin{align*}
\begin{cases}
\displaystyle     1<\dim_{\mathcal{H}}((A\cap B)^2)\leq\frac{57}{56}\\[9pt]
\displaystyle     \dim_{\mathcal{H}}(A)+\frac{(d-2)}{d}\dim_{\mathcal{H}}(B)\geq\frac{21d^2-29d+29}{21d(d-2)}
\end{cases}
{\rm\ \  or\ \ }
\end{align*}
\begin{align*}
\begin{cases}
\displaystyle
    \frac{57}{56}<\dim_{\mathcal{H}}((A\cap B)^2)\leq\frac{5}{4}\\[9pt]
\displaystyle   \left[\frac{(d-2)}{d}\dim_{\mathcal{H}}(B)+\dim_{\mathcal{H}}(A)\right]+\frac{16d-16}{3d(d-2)}\cdot\dim_{\mathcal{H}}(A\cap B)> 1+\frac{10d-4}{3d(d-2)}
\end{cases}
{\rm\ \  or\ \ }
\end{align*}
\begin{align*}
\begin{cases}
\displaystyle
   \frac{5}{4}<\dim_{\mathcal{H}}((A\cap B)^2)\\[9pt]
\displaystyle \dim_{\mathcal{H}}(A)+\frac{(d-2)}{d}\dim_{\mathcal{H}}(B)>1+\frac{2}{d(d-2)},  
\end{cases}
\end{align*}
then the pinned distance set of the product $A^d$ contains an interval, which also improve the threshold in the case of product sets.


\subsection{Proof of Corollary \ref{mainsub1}}
   Let $d\geq3$ and $y_{d-1}, y_d, x_0, x_1\in A$. Consider the sets in $\mathbb{R}^{d-1}$.
    \begin{align*}
        &E:=\left(A\times\cdots \times A\right)\times \triangle^2_{(y_{d-1},y_{d})}(A^2)\subset\mathbb R^{d-2}\times\mathbb R^1;\,\text{and}\\
        &F:=\left(A\times\cdots \times A\right)\times -\triangle^2_{(x_0,x_1)}(A^2)\subset\mathbb R^{d-2}\times\mathbb R^1.
    \end{align*}
    If $\Phi_{d-1}:\mathbb R^{d-1}\times\mathbb R^{d-1}\longrightarrow\mathbb R$ is the parabolic distance,  then there is a probability measure $\mu_F$ on $F$ such that for $\mu_F$-a.e. $y\in F$ which is of the form  $$\left(y_1,\ldots,y_{d-2},-\left(|x_0-a_0|^2+|x_1-a_1|^2\right)\right),$$
    satisfying
     \begin{align*}
        &(1)\,\,\dim_{\mathcal{H}}(\triangle^\Phi_y(E))\geq\beta,\,\,\text{if}\,\dim_{\mathcal{H}}(E)+\frac{d-2}{d}\dim_{\mathcal{H}}(F)>d-2+\beta,\\
          &(2)\,\,\left|\triangle^\Phi_y(E)\right|>0,\,\,\text{if}\,\dim_{\mathcal{H}}(E)+\frac{d-2}{d}\dim_{\mathcal{H}}(F)>d-1,\\
           &(3)\,\,\left(\triangle^\Phi_y(E)\right)^0\neq\emptyset,\,\,\text{if}\,\dim_{\mathcal{H}}(E)+\frac{d-2}{d}\dim_{\mathcal{H}}(F)>d.
    \end{align*}
    However note that $$\triangle^\Phi_y(E)=\triangle^2_{(y_1,\ldots,y_d)}(A^d)+|x_0-a_0|^2+|x_1-a_1|^2,$$
    which says the set $\triangle^2_{(y_1,\ldots,y_d)}(A^d)$ is a translation of $\triangle^\Phi_y(E)$ so that we have
  \begin{align*}
    &(1)\,\,\dim_{\mathcal{H}}(\triangle^\Phi_y(E))=  \dim_{\mathcal{H}}(\triangle^2_{(y_1,\ldots,y_d)}(A^d))=\dim_{\mathcal{H}}(\triangle_{(y_1,\ldots,y_d)}(A^d)),\,\text{and}\\
    &(2)\,\,\left|\triangle^\Phi_y(E)\right|=\left|\triangle^2_{(y_1,\ldots,y_d)}(A^d)\right|\leq C_A\cdot\left|\triangle_{(y_1,\ldots,y_d)}(A^d)\right|,\,\text{and}\\
    &(3)\,\,\left(\triangle^\Phi_y(E)\right)^0\neq\emptyset\implies\left(\triangle_{(y_1,\ldots,y_d)}(A^d)\right)^0\neq\emptyset.
  \end{align*}
Then one has
 \begin{align*}
        &(1)\,\,\dim_{\mathcal{H}}(\triangle_{(y_1,\ldots,y_d)}(A^d))\geq\beta,\,\,\text{if}\,\dim_{\mathcal{H}}(E)+\frac{d-2}{d}\dim_{\mathcal{H}}(F)>d-2+\beta,\\
          &(2)\,\,\left|\triangle_{(y_1,\ldots,y_d)}(A^d)\right|>0,\,\,\text{if}\,\dim_{\mathcal{H}}(E)+\frac{d-2}{d}\dim_{\mathcal{H}}(F)>d-1,\\
           &(3)\,\,\left(\triangle_{(y_1,\ldots,y_d)}(A^d)\right)^0\neq\emptyset,\,\,\text{if}\,\dim_{\mathcal{H}}(E)+\frac{d-2}{d}\dim_{\mathcal{H}}(F)>d.
    \end{align*}
To guarantee the condition $\dim_{\mathcal{H}}(E)+\frac{d-2}{d}\dim_{\mathcal{H}}(F)>d-2+\beta$ holds, we note that
\begin{align}\label{last2}
    &\dim_{\mathcal{H}}(E)+\frac{d-2}{d}\dim_{\mathcal{H}}(F)\notag\\
    &\geq (d-2)\dim_{\mathcal{H}}(A)+\dim_{\mathcal{H}}(\triangle_{(y_{d-1},y_d)}(A^2))\notag\\
    &\qquad+\frac{d-2}{d}\cdot\left((d-2)\dim_{\mathcal{H}}(A)+\dim_{\mathcal{H}}(\triangle_{(x_0,x_1)}(A^2))\right)\notag\\
    &=\frac{(2d-2)\cdot(d-2)}{d}\dim_{\mathcal{H}}(A)+\dim_{\mathcal{H}}(\triangle_{(y_{d-1},y_d)}(A^2))\notag\\
    &\qquad+\frac{d-2}{d}\dim_{\mathcal{H}}(\triangle_{(x_0,x_1)}(A^2)).
\end{align}
To apply Lemma \ref{liu} and Lemma \ref{sm}, we need to assume that the Hausdorff dimension of the set $A$ is at least $\frac{1}{2}$, and hence we have the following three cases: \\[10pt]
\textbf{Case}$(1)$: Assume $1<\dim_{\mathcal{H}}(A^2)\leq\frac{57}{56}$, then (\ref{last2}) reveals that if we choose the point $(y_{d-1},y_d)=(x_0,x_1)$ such that Lemma \ref{sm} holds, then
\begin{align*}
    \dim_{\mathcal{H}}(E)+\frac{d-2}{d}\dim_{\mathcal{H}}(F)> \frac{(2d-2)\cdot(d-2)}{d}\dim_{\mathcal{H}}(A)+\frac{2d-2}{d}\cdot\frac{29}{42}.   
\end{align*}
Suppose that $
\frac{(2d-2)\cdot(d-2)}{d}\dim_{\mathcal{H}}(A)+\frac{2d-2}{d}\cdot\frac{29}{42}\geq d-2+\beta$,
then we have 
\begin{align*}
\dim_{\mathcal{H}}(A)\geq\frac{21d^2-71d+21\beta\cdot d+29}{42(d-1)(d-2)},
\end{align*}
whenever $d\geq3$.\\[10pt] 
\textbf{Case}$(2)$: Assume $\frac{57}{56}<\dim_{\mathcal{H}}(A^2)\leq\frac{5}{4}$, then (\ref{last2}) reveals that if we choose the point $(y_{d-1},y_d)=(x_0,x_1)$ such that Lemma \ref{liu} holds, then
\begin{align*}
   \dim_{\mathcal{H}}(E)+\frac{d-2}{d}\dim_{\mathcal{H}}(F)&\geq \frac{(2d-2)\cdot(d-2)}{d}\dim_{\mathcal{H}}(A)+\frac{2d-2}{d}\cdot\left(\frac{8}{3}\dim_{\mathcal{H}}(A)-\frac{2}{3}\right)\\
&=\frac{(d-1)(6d+4)}{3d}\dim_{\mathcal{H}}(A)-\frac{4(d-1)}{3d}.    
\end{align*}
Suppose that $
\frac{(d-1)(6d+4)}{3d}\dim_{\mathcal{H}}(A)-\frac{4(d-1)}{3d}> d-2+\beta$,
then we have 
\begin{align*}
\dim_{\mathcal{H}}(A)>\frac{3d(d-2+\beta)}{(d-1)(6d+4)}+\frac{2}{3d+2},
\end{align*}
whenever $d\geq3$.\\[10pt] 
\textbf{Case}$(3)$: Assume $\frac{5}{4}<\dim_{\mathcal{H}}(A^2)$, then (\ref{last2}) reveals that if we choose the point $(y_{d-1},y_d)=(x_0,z_0)$ such that Lemma \ref{liu} holds, then
\begin{align*}
   \dim_{\mathcal{H}}(E)+\frac{d-2}{d}\dim_{\mathcal{H}}(F)&\geq \frac{(2d-2)\cdot(d-2)}{d}\dim_{\mathcal{H}}(A)+\frac{2d-2}{d}. 
\end{align*}
Suppose that $
\frac{(2d-2)\cdot(d-2)}{d}\dim_{\mathcal{H}}(A)+\frac{2d-2}{d}> d-2+\beta$,
then we have 
\begin{align*}
\dim_{\mathcal{H}}(A)>\frac{d^2-4d+\beta d+2}{(d-2)(2d-2)},
\end{align*}
whenever $d\geq3$. Therefore, we can conclude that for all $d\geq3$, if
 \begin{align*}
\begin{cases}
\displaystyle     1<\dim_{\mathcal{H}}(A^2)\leq\frac{57}{56}\\[9pt]
\displaystyle     \dim_{\mathcal{H}}(A)\geq\frac{21d^2-71d+21\beta\cdot d+29}{42(d-1)(d-2)}
\end{cases}
{\rm\ \  or\ \ }
\begin{cases}
\displaystyle
    \frac{57}{56}<\dim_{\mathcal{H}}(A^2)\leq\frac{5}{4}\\[9pt]
\displaystyle   \dim_{\mathcal{H}}(A)>\frac{3d(d-2+\beta)}{(d-1)(6d+4)}+\frac{2}{3d+2}
\end{cases}
\end{align*}
or
\begin{align*}
\begin{cases}
\displaystyle
   \frac{5}{4}<\dim_{\mathcal{H}}(A^2)\\[9pt]
\displaystyle \dim_{\mathcal{H}}(A)>\frac{d^2-4d+\beta d+2}{(d-2)(2d-2)},  
\end{cases}
\end{align*}
then the Hausdorff dimension of pinned distance set of the product $A^d$ is no less than $\beta$, which improve the threshold in the case of product sets. Similarly, we have
for all $d\geq3$, if
\begin{align*}
\begin{cases}
\displaystyle     1<\dim_{\mathcal{H}}(A^2)\leq\frac{57}{56}\\[9pt] 
\displaystyle     \dim_{\mathcal{H}}(A)\geq\frac{21d-29}{42(d-2)}
\end{cases}
,\text{or}\,\,
\begin{cases}
\displaystyle    \frac{57}{56}<\dim_{\mathcal{H}}(A^2)\leq\frac{5}{4}\\[9pt]
\displaystyle   \dim_{\mathcal{H}}(A)>\frac{3d+4}{6d+4}
\end{cases}
\end{align*}
or
\begin{align*}
\begin{cases}
  \displaystyle \frac{5}{4}<\dim_{\mathcal{H}}(A^2)\\[9pt]
\displaystyle \dim_{\mathcal{H}}(A)>\frac{1}{2},  
\end{cases}
\end{align*}
then the  pinned distance set of the product $A^d$ has positive Lebesgue measure; and if
\begin{align*}
\begin{cases}
\displaystyle
     1<\dim_{\mathcal{H}}(A^2)\leq\frac{57}{56}\\[9pt]
 \displaystyle    \dim_{\mathcal{H}}(A)\geq\frac{21d^2-29d+29}{42(d-1)(d-2)}
\end{cases}
,\text{or}\,\,
\begin{cases}
\displaystyle    \frac{57}{56}<\dim_{\mathcal{H}}(A^2)\leq\frac{5}{4}\\[9pt]
\displaystyle   \dim_{\mathcal{H}}(A)>\frac{(3d-2)(d+2)}{(6d+4)(d-1)}
\end{cases}
\end{align*}
or
\begin{align*}
\begin{cases}
\displaystyle   \frac{5}{4}<\dim_{\mathcal{H}}(A^2)\\[9pt]
\displaystyle \dim_{\mathcal{H}}(A)>\frac{d^2-2d+2}{(d-2)(2d-2)},  
\end{cases}
\end{align*}
then the pinned distance set of the product $A^d$ contains an interval, which also improve the threshold in the case of product sets.


\subsection{Proof of Corollary \ref{sec}}
  Let $d\geq3$ and $y_{d-1}\in A_{d-1},\, y_d\in A_d$. Consider the sets in $\mathbb{R}^{d-1}$.
    \begin{align*}
        &G:=\left(A_1\times\cdots \times A_{d-2}\right)\times \triangle^2_{(y_{d-1},y_{d})}(A_{d-1}\times A_d)\subset\mathbb R^{d-2}\times\mathbb R^1;\,\text{and}\\
        &F:=\left(A_1\times\cdots \times A_{d-2}\right)\times -\triangle^2_{(y_{d-1},y_{d})}(A_{d-1}\times A_d)\subset\mathbb R^{d-2}\times\mathbb R^1.
    \end{align*}
    If $\Phi_{d-1}:\mathbb R^{d-1}\times\mathbb R^{d-1}\longrightarrow\mathbb R$ is the parabolic distance,  then there is a probability measure $\mu_F$ on $F$ such that for $\mu_F$-a.e. $y\in F$ which is of the form  $$\left(y_1,\ldots,y_{d-2},-\left(|y_{d-1}-a_0|^2+|y_d-a_1|^2\right)\right),$$
    satisfying 
     \begin{align*}
     \left|\triangle^\Phi_y(G)\right|>0,\,\,\text{if}\,\dim_{\mathcal{H}}(G)+\frac{d-2}{d}\dim_{\mathcal{H}}(F)>d-1.
    \end{align*}
    However note that $$\triangle^\Phi_y(G)=\triangle^2_{(y_1,\ldots,y_d)}(A_1\times\cdots\times A_d)+|y_{d-1}-a_0|^2+|y_d-a_1|^2,$$
    which says the set $\triangle^2_{(y_1,\ldots,y_d)}(A_1\times\cdots\times A_d)$ is a translation of $\triangle^\Phi_y(G)$ so that we have
 \begin{align*}
    \left|\triangle^\Phi_y(G)\right|=\left|\triangle^2_{(y_1,\ldots,y_d)}(A_1\times\cdots\times A_d)\right|\leq C_A\cdot\left|\triangle_{(y_1,\ldots,y_d)}(A_1\times\cdots\times A_d)\right|.
  \end{align*}
Then one has
 \begin{align*}
       \left|\triangle^2_{(y_1,\ldots,y_d)}(A_1\times\cdots\times A_d)\right|>0,\,\,\text{if}\,\dim_{\mathcal{H}}(G)+\frac{d-2}{d}\dim_{\mathcal{H}}(F)>d-1,
    \end{align*}
    To guarantee that the Hausdorff dimension of $G$ and $F$ fit the threshold, we may hope that
    \begin{align}\label{aaaa}
        \frac{2d-2}{d}\left(\sum^{d-2}_{j=1}s_j+\dim_{\mathcal{H}}\left(\triangle^2_{(y_{d-1},y_y)}(A_{d-1}\times A_d)\right)\right)>d-1;
    \end{align}
    By lemma \ref{liu}, we have if $s_{d-1}+s_d>\frac{5}{4}$, then $\dim_{\mathcal{H}}(A_{d-1}\times A_d)>\frac{5}{4}$ and hence
    \begin{align}\label{bbbb}
    \dim_{\mathcal{H}}\left(\triangle^2_{(y_{d-1},y_y)}(A_{d-1}\times A_d)\right)=\dim_{\mathcal{H}}\left(\triangle_{(y_{d-1},y_y)}(A_{d-1}\times A_d)\right)\geq 1.
    \end{align}
    Combining (\ref{aaaa}) and (\ref{bbbb}), we have if $\sum^{d-2}_{j=1}s_j>\frac{d}{2}-1$, 
    then $ \left|\triangle^2_{(y_1,\ldots,y_d)}(A_1\times\cdots\times A_d)\right|>0$.

\bigskip

{\bf Acknowledgments:} J. Li is supported by ARC DP 220100285. C.-W. Liang and C.-Y. Shen are supported in part by NSTC through grant 111-2115-M-002-010-MY5. 
C.-W. Liang is also supported by MQ Cotutelle PhD scholarhsip.  

\bigskip

(J Li) Department of Mathematics, Macquarie University, NSW, 2109, Australia.\\ 
{\it E-mail}: \texttt{ji.li@mq.edu.au}\\

(C-W Liang) Department of Mathematics, National Taiwan University, Taiwan.
\\ 
{\it E-mail}: \texttt{d10221001@ntu.edu.tw}\\
 
(C-Y Shen) Department of Mathematics, National Taiwan University, Taiwan. \\ {\it E-mail}: \texttt{cyshen@math.ntu.edu.tw}

\end{document}